\documentclass[12pt]{article}
\usepackage[cp1251]{inputenc}
\usepackage[russian]{babel}
\usepackage{amssymb,amsthm,amsmath}

\theoremstyle{plain}

\newtheorem{Theorem}{Теорема}
\newtheorem{Statement}{Утверждение}

\newcommand{\intl}{\mathop{\int}\limits}
\newcommand{\liml}{\mathop{\lim}\limits}
\newcommand{\suml}{\mathop{\sum}\limits}

\renewcommand{\le}{\leqslant}
\renewcommand{\ge}{\geqslant}

\def\la{\lambda}

\begin{document}

{\bf УДК 517.984}

{\large {\bf \centerline{О скорости равносходимости разложений в
ряды}} {\bf \centerline{по тригонометрической системе и по
собственным функциям}} {\bf \centerline{оператора Штурма--Лиувилля
с потенциалом --- распределением }} }

\bigskip

{\centerline{Садовничая~И.~В.}}

\medskip
В настоящей работе изучается оператор Штурма--Лиувилля
\begin{equation}\label{main}
Ly=l(y)=-\dfrac{d^2y}{dx^2}+q(x)y,
\end{equation}
в пространстве \(L_2[0,\pi]\) с граничными условиями Дирихле
$y(0)=y(\pi)=0$. Предполагается, что потенциал  \(q(x)=u'(x),\)
\(u\in W_2^{\theta}[0,\pi]\), \(0<\theta<1/2\). Производная здесь
понимается в смысле распределений, а через $W_2^{\theta}[0,\pi]$
мы обозначаем соболевские пространства с дробным показателем (см.,
например, \cite[п. 4.2]{Tr}). Операторы такого вида были
определены в работе \cite{SSh1}. В работах \cite{SSh1} --
\cite{SSh2} было доказано, что оператор \(L\) фредгольмов с
индексами (0,0) (в случае вещественного потенциала --
самосопряжен), полуограничен, имеет чисто дискретный спектр.
Некоторые результаты об операторе \(L\), полученные в \cite{SSh1}
-- \cite{SSh3} и необходимые нам в данной работе, будут приведены
ниже.

В статье рассматривается вопрос о равномерной на всем отрезке
\([0,\pi]\) равносходимости разложения функции \(f\) в ряд по
системе собственных и присоединенных функций опрератора \(L\) с ее
разложением в ряд Фурье по системе синусов. Эта задача хорошо
известна в классической теории операторов Штурма--Лиувилля (в
случае, когда потенциал локально суммируем). В монографии
В.~А.~Марченко \cite[\S 3, гл. 1]{Mar} была доказана равномерная
равносходимость в случае, если $f\in L_2[0,\pi]$, а $q$ ---
комплекснозначная суммируемая функция. В 1985 году В.~А.~Ильин
получил результаты в случае, когда $f\in L_1[0,\pi]$,  $q$ ---
комплекснозначная суммируемая функция (см. \cite{Il}).
В.~А.~Винокуров и В.~А.~Садовничий в \cite{VS} доказали теорему о
равносходимости для случая операторов с потенциалом ---
производной функции ограниченной вариации, при этом \(f\in
L_1[0,\pi]\).

Вопрос о скорости равносходимости для классических потенциалов
изучался в статье А.~М.~Гомилко и Г.~В.~Радзиевского \cite{R}.
Здесь мы покажем, что, если первообразная \(u\) от потенциала ---
комплекснозначная функция из пространства \(W_2^{\theta}[0,\pi]\),
\(0<\theta<1/2\), то для любой функции \(f\) из пространства
\(L_2[0,\pi]\) можно оценить скорость равносходимости равномерно
по шару \(u\in B_{\theta,R}=\{v\in
W_2^\theta[0,\pi]|\|v\|_{W_2^\theta}\le R\}\). Оценки такого вида
являются новыми даже для классического случая $q\in L_2[0,\pi]$.

Случай \(\theta=0\) является особым и требует отдельного изучения.
Некоторые результаты о равносходимости в случае вещественного
сингулярного потенциала были анонсированы автором в работе
\cite{Sad}, но их подробное доказательство опубликовано не было.
Случай комплексного потенциала $u\in L_2[0,\pi]$ автор планирует
рассмотреть в другой работе.

\newpage
{\bf \S 1. Предварительные результаты.}

\medskip

 Нам понадобятся некоторые предварительные сведения об
операторе \eqref{main}. Обозначим через $\omega(x,\la)$ решение
дифференциального уравнения $-\omega''+q\omega=\la\omega$ с
начальными условиями $\omega(0,\la)=0$, $\omega^{[1]}(0,\la)=1$
(здесь $\omega^{[1]}=\omega'-u\omega$ --- первая
квазипроизводная). Ясно, что нули целой функции $\omega(\pi,\la)$
совпадают с собственными значениями оператора $L$. Алгебраической
кратностью собственного значения $\la_0$ мы будем называть
кратность нуля $\la_0$ функции $\omega(\pi,\la)$. В силу теоремы
существования и единственности, геометрическая кратность каждого
собственного значения равна $1$. Пусть $\la_0$ есть собственное
значение алгебраической кратности $p\ge 2$. Заметим, что функции
$\omega^{(j)}_{\la}(x,\la)$, $j=1,2,\dots,p-1$ удовлетворяют
дифференциальным уравнениям
$l(\omega^{(j)}_{\la})=\la\omega^{(j)}_{\la}+\omega^{(j-1)}_\la$,
причем $\omega^{(j)}_{\la}(0,\la)\equiv0$. Кроме того,
$\omega^{(j)}_\la(\pi,\la_0)=0$. Таким образом, функции
$\omega^{(j)}_\la(x,\la_0)$, $j=0,1,\dots,p-1$ образуют цепочку из
собственной и присоединенных функций, отвечающую собственному
значению $\la_0$. Легко видеть, что функции этой системы линейно
независимы, а значит порождают подпространство размерности $p$.
Обозначим собственные значения $\{\la_n\}_{n=1}^\infty$ и
занумеруем их в порядке возрастания модуля
$|\la_1|\le|\la_2|\le\dots$ с учетом алгебраической кратности. В
случае совпадения модулей нумерацию будем вести по возрастанию
аргумента, значения которого выбираются из полуинтервала
$(-\pi,\pi]$. Под системой собственных и присоединенных функций
$\{y_n(x)\}_{n=1}^\infty$ мы будем понимать систему, полученную
нормировкой $\|y_n\|_{L_2}=1$ системы
$\bigcup\limits_{k=1}^{\infty}\{\omega^{(j)}_\la(x,\la_k)\}_{j=0}^{p_k}$
(здесь $\lambda_k$
--- все различные собственные значения, $p_k$
--- алгебраическая кратность $\la_k$). Отметим, что в случае
вещественного потенциала все собственные значения являются
простыми.

\begin{Statement} (см. теорему 2.9 работы \cite{SSh3})\label{st:5}
Пусть $u\in W_2^\theta[0,\pi]$ при некотором $\theta\in[0,1/2)$.
Тогда система $\{y_n(x)\}_{n=1}^\infty$ собственных и
присоединенных функций оператора \(L\) образует базис Рисса в
пространстве $L_2[0,\pi]$.
\end{Statement}
Таким образом, существует единственная биортогональная система
$\{w_n(x)\}_{n=1}^\infty$ (т.е. система, для которой выполнены
равенства $(y_n(x),w_m(x))=\delta_n^m$). \footnote{всюду в
дальнейшем символом \((\cdot,\cdot)\) будет обозначаться скалярное
произведение в пространстве \(L_2[0,\pi]\); биортогональность
понимается также в \(L_2[0,\pi]\)} В нашем случае систему
$\{w_n(x)\}_{n=1}^\infty$ можно выписать явно в виде конечных
линейных комбинаций собственных и присоединенных функций оператора
$L^*=-d^2/dx^2+\overline {q(x)}$ (см. \cite{S}). В частности,
$w_n(x)=\overline{y_n(x)}/(y_n(x),\overline{y_n(x)})$, если
собственное значение $\lambda_n$ имеет алгебраическую кратность 1.

Обозначим через $l_2^\theta$ пространство последовательностей
$\{\alpha_n\}_{n=1}^\infty$, для которых
$\left(\suml_{n=1}^\infty|\alpha_n|^2n^{2\theta}\right)^{1/2}=\|\{\alpha_n\}\|_{l_2^\theta}<
\infty$. Напомним, что для любой функции $u\in W_2^\theta[0,\pi]$,
$0<\theta<1/2$, выполнено:
$C_1\|\{u_n\}\|_{l_2^\theta}\le\|u\|_{W_2^\theta}\le
C_2\|\{u_n\}\|_{l_2^\theta}$, где $u_n=(2/\pi)(u(x),\sin nx)$.
\begin{Statement}\label{st:7}
Пусть $u\in W_2^\theta[0,\pi]$, $0<\theta<1/2$, а
$\|u\|_{W_2^\theta}\le R$. Тогда, начиная с некоторого номера
$N=N_{\theta,R}$ (нижние индексы здесь и в дальнейшем указывают на
то, что выбор номера $N$ зависит только от $\theta$ и $R$ и не
зависит от других характеристик потенциала), все собственные
значения оператора $L$ просты, а для функций $y_n$ и $w_n$
справедливы асимптотические равенства
\begin{equation}\label{efas}
\begin{array}{c}
y_n(x)=\sqrt{\frac2\pi}\sin nx+\varphi_n(x),\quad
w_n(x)=\sqrt{\frac2\pi}\sin nx+\psi_n(x),\\
y'_n(x)=n\left(\sqrt{\frac{2}{\pi}}\cos nx+\eta_n(x)\right)+
u(x)\left(\sqrt{\frac{2}{\pi}}\sin nx+\varphi_n(x)\right)\quad
n=N,N+1,\dots\ ,
\end{array}
\end{equation}
причем последовательность
$\{\gamma_n\}_{n=N}^\infty=\{\|\varphi_n(x)\|_{C}+
\|\psi_n(x)\|_{C}+\|\eta_n(x)\|_{C}\}_{n=N}^\infty\in l^\theta_2$
и ее норма в этом пространстве ограничена постоянной, зависящей
только от $\theta$ и $R$. Кроме этого, при $n\ge N$
\begin{equation} \label{efas1} \psi_n(x)=\alpha_n \sin
nx+\beta_n\cos nx-\intl_0^xu(t)\sin n(x-2t)dt+\psi_n^1(x),
\end{equation} где
норма
 последовательности
$\{|\alpha_n|+|\beta_n|\}_{n=N}^\infty$ в пространстве
$l_2^{\theta}$, а также норма последовательности
\(\{\|\psi_n^1(x)\|_{C}\}_{n=N}^\infty\) в $ l_1$ ограничены
величиной $C_{\theta, R}$.
\end{Statement}

Доказательство этого утверждения (с оценками остатков, зависящими
от потенциала) было приведено в работе \cite{SSh2} (теорема 3.13).
Это доказательство опиралось на асимптотические формулы для
собственных значений оператора $L$, полученные в той же статье.
Однако позже, в работе \cite{SSh3}, были выведены асимптотики
собственных значений с оценками остатков, равномерными по шару
\(u\in B_{\theta,R}\). С учетом этих результатов, рассуждениями,
полностью аналогичными доказательству теоремы 3.13 работы
\cite{SSh2}, можно получить формулы \eqref{efas}--\eqref{efas1}.
\begin{Statement}\label{st:8}
Пусть $u\in W_2^\theta[0,\pi]$, где $0<\theta<1/2$, и
$\|u\|_{W_2^\theta}\le R$. Тогда существует такой номер
$N_{\theta, R}$, что для любого $n\ge N_{\theta, R}$ оператор
проектирования $P_n(u)$, определенный по правилу
$P_nf=\suml_{k=1}^n(f,w_n)y_n$ и действующий из пространства
$L_2[0,\pi]$ в пространство $W_2^1[0,\pi]$, непрерывно зависит от
параметра $u$. А именно, $\|P_n(u)-P_n(u_0)\|_{L_2\to W_2^1}\to0$,
если $\|u-u_0\|_{W_2^\theta}\to0$.
\end{Statement}

Это утверждение вытекает из результатов работы \cite{SSh2}
(теорема 1.9) и классических результатов о полунепрерывности
изолированных частей спектра (см. \cite[теоремы IV.2.23 и
IV.3.16]{Ka}). Более подробное доказательство этого факта
приведено в \cite{S}.

Отметим, что для произвольного натурального $N$ это утверждение,
вообще говоря, неверно. Действительно, пусть $R$ достаточно
велико, функция $u$ комплекснозначна, а $\la_n(u)$ и
$\la_{n+1}(u)$
--- пара простых собственных значений. При изменении функции $u$ в
шаре $B_{\theta, R}$ эти собственные значения могут начать
сближаться и затем, при некотором $u_0$, столкнуться, образовав
клетку. Можно проверить, что в этом случае число
$(y_n,\overline{y_n})\to 0$ при $u\to u_0$. Мы уже отмечали, что
для простых собственных значений вектор $w_n$ биортогональной
системы имеет вид
$w_n(x)=\overline{y_n(x)}/(y_n(x),\overline{y_n(x)})$. Таким
образом, норма вектора $w_n$, а значит, и одномерного проектора
$P_nf=(f,w_n)y_n$, неограниченно растет при $u\to u_0$.

{\bf \S 2. Основная теорема.}
\begin{Theorem}
Пусть \(R>0\). Рассмотрим оператор \eqref{main}, действующий в
пространстве \(L_2[0,\pi]\), с граничными условиями Дирихле,
потенциал которого удовлетворяет следующим условиям:\,
 \(q(x)=u'(x)\), где комплекснозначная функция \(u\in W_2^\theta [0,\pi]\), \(0<\theta <
1/2 \), причем \(\|u\|_{W_2^\theta}\le R\). Пусть \(\{
y_n(x)\}_{n=1}^\infty\) -- нормированная система собственных и
присоединенных функций оператора \(L\), \(\{
w_n(x)\}_{n=1}^\infty\) -- биортогональная к ней система. Для
произвольной функции \(f\in L_2 [0,\pi]\) обозначим
\(c_n=(f(x),w_n(x))\), \(c_{n,0}=\sqrt{2/\pi}(f(x),\sin{nx})\).
Тогда существует натуральное число \(M=M_{\theta,R}\) такое, что
для любого \(m\ge M\)

\begin{equation}
\label{1}
\left\|\sum_{n=1}^{m}{c_n
y_n(x)}-\sum_{n=1}^{m}\sqrt{\frac{2}{\pi}}{c_{n,0}\sin
nx}\right\|_{C}\le C_{\theta,
R,\varepsilon}\left(\sqrt{\suml_{n\ge\sqrt{m}}|c_{n,0}|^2}+\frac{\|f\|_{L_2}}
{m^{\theta/2-\varepsilon}}\right).
\end{equation}

Здесь \(\varepsilon\)-- сколь угодно малое положительное число.
\end{Theorem}

{\bf Доказательство теоремы 1.}\quad Рассмотрим операторы
\(B_{m,N}: L_2[0,\pi]\to C[0,\pi]\), действующие по правилу
$$
B_{m,N}f(x):=\sum_{n=N}^{m}{c_ny_n(x)}-\sum_{n=N}^{m}\sqrt{\frac{2}{\pi}}{c_{n,0}\sin
nx},
$$ где \(f\in L_2[0,\pi].\) При $N=1$ будем опускать второй
индекс: $B_{m,1}=:B_m$.

Пусть теперь $N=N_{\theta,R}$ -- наибольшее из двух натуральных
чисел, существование которых постулируется в утверждениях 2 и 3.

\medskip
{\bf  Шаг 1 (оценка нормы оператора $B_{m,N}$).}
\medskip

{\it Пусть \(R>0,\) \(\|u\|_{W_2^{\theta}}\le R\), \(0<\theta
<1/2\).  Тогда
\begin{equation} \label{B}
\|B_{m,N}(u)\|_{L_2\to C}\le C_{\theta, R}.
\end{equation}}

Этот шаг является ключевым для доказательства теоремы 1.
Естественно, он наиболее сложен.

Очевидно, что для любой функции $f\in L_2[0,\pi]$
\begin{equation}
\label{3}
B_{m,N}f(x)=\sum_{n=N}^{m}\sqrt{\frac{2}{\pi}}(f(t),\psi_n(t))\sin
nx+ \sum_{n=N}^{m}\sqrt{\frac{2}{\pi}}(f(t),\sin
nt)\varphi_n(x)+\sum_{n=N}^{m}(f(t),\psi_n(t))\varphi_n(x).
\end{equation}

Оценим каждое из слагаемых в правой части соотношения \eqref{3} по
отдельности (наиболее тяжелой здесь будет оценка первого
слагаемого). В силу асимптотических формул \eqref{efas1}
$$
(f(t),\psi_n(t))=\overline{\alpha_n}(f(t),\sin
nt)+\overline{\beta_n}(f(t),\cos nt)-
$$
\begin{equation}
\label{4} -\intl_0^{\pi}f(t)\intl_0^{t}\overline{u(s)}\sin
n(t-2s)dsdt+\intl_0^{\pi}f(t)\overline{\psi_n^1(t)}dt.
\end{equation}
 Так как $\|\{\alpha_n\}\|_{l_2^{\theta}}\le C_{\theta, R}$, то
 $$
\left
\|\sum_{n=N}^{m}\sqrt{\frac{2}{\pi}}\overline{\alpha_n}(f(t), \sin
nt)\sin nx\right\|_{C}\le C_{\theta, R}\|f\|_{L_2}.
 $$
 Аналогичные рассуждения справедливы и для второго слагаемого в
 \eqref{4}. Наиболее сложными для оценки является третье слагаемое.
 Применив формулу разности косинусов, получим
 $$
\left\|\sum_{n=N}^{m}\sin
nx\intl_0^{\pi}f(t)\intl_0^{t}\overline{u(s)}\sin
n(t-2s)dsdt\right\|_{C}=\frac{1}{2}\left\|\intl_0^{\pi}f(t)\intl_0^{t}\overline{u(s)}
\sum_{n=N}^{m}(\cos n(t-2s+x)- \right.$$ $$-\cos
n(t-2s-x))dsdt\Bigg\|_{C}\le\frac{1}{2}\left\|\intl_0^{\pi}f(t)\intl_0^{t}
\overline{u(s)}(D_m(t-2s+x)-D_m(t-2s-x))dsdt\right\|_{C}+
$$
$$
+\frac{1}{2}\left\|\intl_0^{\pi}f(t)\intl_0^{t}
\overline{u(s)}(D_{N-1}(t-2s+x)-D_{N-1}(t-2s-x))dsdt\right\|_{C},
$$
где $D_m(\xi)=1/2+\suml_{n=1}^{m}\cos n\xi$ -- ядро Дирихле.
Заметим, что $N=N_{\theta,R}$ и, следовательно, последнее
слагаемое не превосходит $ C_{\theta,R}\|f\|_{L_2}$. Разобьем
первое слагаемое на сумму двух и оценим один из интегралов (второй
можно оценить аналогично). $$ \left\|\intl_0^{\pi}f(t)\intl_0^{t}
\overline{u(s)}D_m(t-2s+x)dsdt\right\|_{C}\le\left\|\intl_0^{\pi}\intl_{x-t}^{x+t}f(t)
\overline{u}\left(\frac{t+x}{2}\right)D_m(\xi)d\xi
dt\right\|_{C}+$$
$$+\left\|\intl_0^{\pi}\intl_{x-t}^{x+t}f(t)
\left(\overline{u}\left(\frac{t-\xi +x}{2}\right)\right.\right.
\left.-\overline{u}\left(\frac{t+x}{2}\right)\right)D_m(\xi)d\xi
dt\Bigg\|_{C}=I_1+I_2
$$
Заметим, что $I_1\le C_{\theta,R}\|f\|_{L_2}, $ поскольку
$\left|\intl_{x-t}^{x+t}D_m(\xi)d\xi\right|\le C$ (см., например,
\cite[гл. 1 \S 35]{Ba}).

Далее, \begin{equation} \label{4'}I_2\le
\left(\intl_0^{\pi}|f(t)|^2dt\right)^{1/2}\cdot
\left\|\intl_0^{\pi}\left|\intl_{x-t}^{x+t}
\left(u\left(\frac{t-\xi
+x}{2}\right)-u\left(\frac{t+x}{2}\right)\right)D_m(\xi)d\xi\right|^{2}
dt\right\|^{1/2}_{C}. \end{equation}

Рассмотрим внутренний интеграл во втором сомножителе правой части
\eqref{4'}. Так как $|D_m(\xi)|\le C/|\xi|$ (\cite[гл. 1 \S
32]{Ba}), то
$$ \left|\intl_{x-t}^{x+t}
\left(u\left(\frac{t-\xi
+x}{2}\right)-u\left(\frac{t+x}{2}\right)\right)D_m(\xi)d\xi\right|
\le C\intl_{x-t}^{x+t}\frac{1}{|\xi|}\left|u\left(\frac{t-\xi
+x}{2}\right)-u\left(\frac{t+x}{2}\right)\right|d\xi\le
$$
$$
\le
C\left(\intl_{x-t}^{x+t}\frac{1}{|\xi|^{1-2\theta}}d\xi\right)^{1/2}
\left(\intl_{x-t}^{x+t}\frac{1}{|\xi|^{1+2\theta}}\left|u\left(\frac{t-\xi
+x}{2}\right)-u\left(\frac{t+x}{2}\right)\right|^2d\xi\right)^{1/2}.
$$

Отсюда и из \eqref{4'} получаем, что \begin{equation}\label{4''}
I_1+I_2\le
C_{\theta,R}\|f\|_{L_2}\left(1+\left(\intl_{0}^{\pi}\intl_{-\pi}^{\pi}
\frac{|u(y-\eta)-u(y)|^2}{|\eta|^{1+2\theta}}d\eta
dy\right)^{1/2}\right).
\end{equation}
Мы сделали замену переменной в последнем интеграле и расширили
границы интегрирования с учетом того, что под знаком интеграла
стоит неотрицательная функция.\footnote{всюду в оценках интегралов
при необходимости полагаем, что все функции периодически
продолжены за отрезок $[0,\pi]$.}

Осталось оценить двойной интеграл в последнем выражении. Поскольку
сумма $\|u\|_{L_2}+\sqrt{\intl_{|h|\le
\delta}\frac{\|u(x+h)-u(x)\|^2_{L_2}}{|h|^{2\theta+1}}dh}$
определяет эквивалентную норму в пространстве $W_2^\theta[0,\pi]$
(см. \cite[п. 4.4.2]{Tr}), то сразу получаем необходимую оценку.
Однако для удобства читателей приведем здесь также простое
доказательство этого факта. Обозначим
$\Psi(\eta)=\intl_0^\pi|u(y-\eta)-u(y)|^2dy$. Так как
$$u(y-\eta)-u(y)=\suml_{n=1}^{\infty}u_n(\sin n(y-\eta)-\sin
ny)=\suml_{n=1}^{\infty}u_n(\cos n\eta-1)\sin
ny-\suml_{n=1}^{\infty}u_n\sin n\eta\cos ny,$$ где
$u_n=(2/\pi)(u(y),\sin ny)$, то $$\Psi(\eta)\le
\suml_{n=1}^{\infty}|u_n|^2(\cos
n\eta-1)^2+\suml_{n=1}^{\infty}|u_n|^2\sin^2
n\eta=4\suml_{n=1}^{\infty}|u_n|^2\sin^2 (n\eta/2).$$ Тогда
$$
\intl_{0}^{\pi}\intl_{-\pi}^{\pi}
\frac{|u(y-\eta)-u(y)|^2}{|\eta|^{1+2\theta}}d\eta
dy=\intl_{-\pi}^\pi\frac{\Psi(\eta)}{|\eta|^{1+2\theta}}d\eta=4\suml_{n=1}^\infty|u_n|^2
\intl_{-\pi}^\pi\frac{\sin^2(n\eta/2)}{|\eta|^{1+2\theta}}d\eta\le
C\suml_{n=1}^\infty|u_n|^2n^{2\theta}.
$$
Поскольку функция $u\in W_2^\theta[0,\pi]$, то последний ряд
сходится и его сумма ограничена константой, зависящей только от
$\theta$ и $R$. Значит, с учетом неравенства \eqref{4''}, можем
получить оценку третьего слагаемого в \eqref{4}:
$$
\left\|\sum_{n=N}^{m}\sin
nx\intl_0^{\pi}f(t)\intl_0^{t}\overline{u(s)}\sin
n(t-2s)dsdt\right\|_{C}\le C_{\theta,R}\|f\|_{L_2}.
$$

Оценка последнего слагаемого в \eqref{4} вытекает непосредственно
из асимптотики \eqref{efas1}:
$$
\sum_{n=N}^{m}\left|\intl_0^{\pi}f(t)\overline{\psi_n^1(t)}dt\right|\le
\sum_{n=N}^{m}\|\psi_n^1(x)\|_{C}\intl_0^{\pi}|f(t)|dt\le
C_{\theta, R}\|f\|_{L_2}.
$$

 Таким образом,
 \begin{equation}
 \label{5}
 \left\|\sum_{n=N}^{m}\sqrt{\frac{2}{\pi}}(f(t),\psi_n(t))\sin
 nx\right\|_{C}\le C_{\theta, R}\|f\|_{L_2}.
\end{equation}

Перейдем ко второму члену представления \eqref{3}. В силу
асимптотических формул \eqref{efas}:
\begin{equation} \label{6}
\left\|\sum_{n=N}^{m}\sqrt{\frac{2}{\pi}}(f(t),\sin
nt)\varphi_n(x)\right\|_{C}\le
\sum_{n=N}^{m}(|f_n|\cdot\|\varphi_n(x)\|_{C})\le C_{\theta,
R}\|f\|_{L_2},
\end{equation}
где $f_n=\sqrt{2/\pi}(f(x),\sin nx)$.

Наконец,
\begin{equation}
\label{7}
\left\|\sum_{n=N}^{m}(f(t),\psi_n(t))\varphi_n(x)\right\|_{C}\le
\|f\|_{L_2}\sum_{n=N}^{m}(\|\psi_n(t)\|_{L_2}\cdot\|\varphi_n(x)\|_{C})\le
C_{\theta, R}\|f\|_{L_2}.
\end{equation}
Из неравенств \eqref{5}--\eqref{7} вытекает оценка \eqref{B}.

Шаг 1 завершен.

\medskip
{\bf  Шаг 2 (оценка нормы оператора $B_{m}$).}
\medskip

{\it Пусть \(R>0,\) \(\|u\|_{W_2^{\theta}}\le R\), \(0<\theta
<1/2\).  Тогда
\begin{equation} \label{B'}
\|B_{m}(u)\|_{L_2\to C}\le C_{\theta, R}.
\end{equation}}

Очевидно, что оператор $B_m$ представляется в виде суммы:
\(B_m=B_{N-1}+B_{m,N}=P_{N-1}(u)-P_{N-1}(0)+B_{m,N}\) (здесь
$P_n(u)$ --- операторы проектирования, введенные в утверждении 3).
Поскольку выбор числа \(N\) зависит только от \(\theta\) и \(R\),
то $\|P_{N-1}(0)\|_{L_2\to C}\le C_{\theta,R}$. Остается оценить
$\|P_{N-1}(u)\|_{L_2\to C}$. Из утверждения 3 вытекает, что норма
оператора проектирования \(P_{N-1}\), определенного по правилу
\(P_{N-1}f=\suml_{n=1}^{N-1}(f,w_n)y_n\) и действующего из
пространства \(L_2[0,\pi]\) в пространство \(W_2^1[0,\pi]\),
ограничена константой, зависящей только от \(\theta\) и \(R\).
Действительно, пусть \(u\in W_2^\theta [0,\pi]\), где \(0<\theta <
1/2 \), причем \(\|u\|_{W_2^\theta}\leq R\). Тогда \(u\in
W_2^{\theta/2} [0,\pi]\) и \(\|u\|_{W_2^{\theta/2}}\leq R\). Так
как вложение шара радиуса \(R\) пространства \(W_2^\theta
[0,\pi]\)  в пространство \(W_2^{\theta/2} [0,\pi]\) компактно, то
непрерывная функция \(\|P_{N-1}(u)\|_{L_2\to W_2^1}\) будет
достигать на нем своих точных граней. В силу теоремы вложения
Соболева: $W_2^{1}[0,\pi]\hookrightarrow C[0,\pi]$ (см., например,
\cite[п. 4.6.2]{Tr}) имеем, что \(\|P_{N-1}(u)\|_{L_2\to C}\le
\|P_{N-1}(u)\|_{L_2\to W_2^1}\le C_{\theta,R}\).

Шаг 2 завершен.

\medskip
{\bf Шаг 3 (доказательство равносходимости).}

\medskip

{\it Для любой функции \(f\in L_2[0,\pi]\) выполнено:
\begin{equation}\label{B1} \lim_{m\to\infty} \|B_mf\|_C=0. \end{equation}}

Рассмотрим действие оператора \(B_m\) на собственные и
присоединенные функции оператора \(L\):
$$
B_my_k(x)=\suml_{n=1}^{m}(y_k(x),w_n(x))y_n(x)-\frac{2}{\pi}\suml_{n=1}
^{m}(y_k(x),\sin nx)\sin nx.
$$
Первое слагаемое в правой части последнего соотношения равно 0 при
\(m<k\) и равно \(y_k(x)\) при \(m\ge k\). Второе слагаемое
представляет собой частичную сумму ряда Фурье функции \(y_k\). Так
как все функции \(y_k\in W_2^1 [0,\pi]\), то ряд Фурье функции
\(y_k\) сходится к ней равномерно на отрезке \([0,\pi]\), и мы
получаем, что $\liml_{m\to\infty}\|B_my_k\|_C=0. $

Осталось заметить, что, в силу полноты системы \(\{y_k(x)\}\) (см.
утверждение 1), из непрерывности оператора \(B_m\) следует
предельное соотношение \eqref{B1}.

Шаг 3 завершен.

\medskip

Перейдем к доказательству утверждения о скорости равносходимости.
Для любого $k\ge N_{\theta, R}$ обозначим
\(g_k(x)=\suml_{n=1}^{k}c_ny_n(x)\) (напомним, что
\(c_n=(f(x),w_n(x))\)). Очевидно, что для любой функции \(f\in
L_2[0,\pi]\) и для любого натурального \(m\) выполнено:
\begin{equation}
\label{9} \|B_m f\|_{C}\le\|B_m(f-g_k)\|_{C}+\|B_mg_k\|_{C}.
\end{equation}

\medskip
{\bf Шаг 4 (оценка нормы $B_m(f-g_k)$ в пространстве $C[0,\pi]$).}

{\it Пусть \(g_k(x)=\suml_{n=1}^{k}c_ny_n(x)\), $k\ge N_{\theta,
R}$. Тогда для любого $m\in\mathbb{N}$
\begin{equation} \label{10} \|B_m(f-g_k) \|_{C}\le
C_{\theta,R}\left(\left(\suml_{n=k+1}^{\infty}|c_{n,0}|^2\right)^{1/2}+\frac{\|f\|_{L_2}}{k^{\theta}}\right),
\end{equation}где $c_{n,0}=\sqrt{2/\pi}(f(x),\sin nx)$.}

С учетом асимптотических формул \eqref{efas} получаем:
$$
\|f(x)-g_k(x)\|_{L_2}\le
\left\|\suml_{n=k+1}^{\infty}\frac{2}{\pi}(f(x),\sin nx)\sin
nx\right\|_{L_2}+
\left\|\suml_{n=k+1}^{\infty}\sqrt{\frac{2}{\pi}}(f(x),
\psi_n(x))\sin nx\right\|_{L_2}+$$
$$+\left\|\suml_{n=k+1}^{\infty}\sqrt{\frac{2}{\pi}}(f(x),\sin
nx)\varphi_n(x)\right\|_{L_2}+ \left\|\suml_{n=k+1}^{\infty}(f(x),
\psi_n(x))\varphi_n(x)\right\|_{L_2}\le$$
$$\le\left(\left(\suml_{n=k+1}^{\infty}|c_{n,0}|^2\right)^{1/2}
+\left(\suml_{n=k+1}^{\infty}|(f(x),\psi_n(x))|^2\right)^{1/2}\right)
\left(1+\frac{C_{\theta,R}}{k^\theta}\right)\le
$$
$$\le
\left(\left(\suml_{n=k+1}^{\infty}|c_{n,0}|^2\right)^{1/2}+
\frac{\|f\|_{L_2}}{k^{\theta}}\right)\left(1+\frac{C_{\theta,R}}{k^\theta}\right).
$$

Поскольку $\|B_m\|_{L_2\to C}\le C_{\theta,R}$ (см. \eqref{B'}),
то из доказанного неравенства немедленно вытекает оценка
\eqref{10}.

Шаг 4 завершен.

Перейдем к оценке второго слагаемого в \eqref{9}. Пусть $m>k$.
Обозначим через $S_m$ оператор, действующий из пространства
$W_2^1[0,\pi]$ в пространство $C[0,\pi]$ по правилу:
$S_mh(x)=2/\pi\suml_{n=m+1}^{\infty}(h(t),\sin nt)\sin nx$.
Заметим, что
$$B_mg_k(x)=g_k(x)-\frac{2}{\pi}\suml_{n=1}^{m}(g_k(t),\sin nt)\sin
nx=S_mg_k(x)$$ (поскольку все собственные и присоединенные функции
оператора $L$ принадлежат пространству $W_2^1[0,\pi]$, то действие
оператора $S_m$ на них корректно определено). Тогда
\begin{equation}
\label{11} \|B_{m} g_k\|_{C}\le
\|S_m(g_k-g_{N})\|_{C}+\|S_mg_{N}\|_{C}\end{equation}

{\bf Шаг 5 (оценка нормы $S_m(g_k-g_{N})$ в пространстве
$C[0,\pi]$).}

{\it Пусть $k,m$ -- натуральные числа, $m>k\ge N_{\theta, R}$.
Тогда для любой функции $f$ из пространства $L_2[0,\pi]$ и любого
числа $\varepsilon>0$ справедлива оценка}
\begin{equation}
\label{S}  \|S_m(g_k-g_{N})\|_{C}\le C_{\theta,
R,\varepsilon}\|f\|_{L_2}m^{\varepsilon-1/2}k^{1-\theta}.
\end{equation}
Заметим, что левую часть неравенства \eqref{S} можно представить в
виде
$$\|S_m(g_k(x)-g_{N}(x))\|_{C}=\left\|\suml_{n=N+1}^{k}c_nS_my_n(x)\right\|_{C}=
\left\|\suml_{n=N+1}^{k}c_nS_m\varphi_n(x)\right\|_{C},$$ так как
$y_n(x)=\sqrt{2/\pi}\sin nx+\varphi_n(x)$, где \(\varphi_n(x)\)
определены в \eqref{efas}. Далее,
$$
\left\|\suml_{n=N+1}^{k}c_nS_m\varphi_n(x)\right\|_{C}
\le\left(\suml_{n=N+1}^{k}|c_n|^2\right)^{1/2}
\left(\suml_{n=N+1}^{k}\|S_m\varphi_n(x)\|_C^2\right)^{1/2}\le
$$ $$\le
C_{\theta,
R}\|f\|_{L_2}\left(\suml_{n=N+1}^{k}\|S_m\varphi_n(x)\|_{C}^2\right)^{1/2}.
$$
При выводе последнего неравенства мы учли тот факт, что
$$\suml_{n=N+1}^{k}|c_n|^2\le \suml_{n=N+1}^{\infty}|c_n|^2\le
2\left(\frac{2}{\pi}\suml_{n=N+1}^{\infty}|(f(x),\sin nx)|^2
+\suml_{n=N+1}^{\infty}|(f(x),\psi_n(x))|^2\right)\le
C_{\theta,R}\|f\|_{L_2}^2. $$

Оценим  $\|S_m\varphi_n(x)\|_C$, воспользовавшись теоремой
вложения Соболева: пространство
$W_2^{\varepsilon+1/2}[0,\pi]\hookrightarrow C[0,\pi]$ (\cite[п.
4.6.2]{Tr}). Получим, что для любого \(\varepsilon>0\),
$$
\|S_m\varphi_n(x)\|_{C}\le
C_{\varepsilon}\left(\suml_{j=m+1}^{\infty}j^{2\varepsilon+1}|(\varphi_n(x),\sin
jx)|^2\right)^{1/2}=
$$
\begin{equation}
\label{12}
=C_{\varepsilon}\left(\suml_{j=m+1}^{\infty}j^{2\varepsilon-1}|(\varphi'_n(x),\cos
jx)|^2\right)^{1/2}\le
C_{\varepsilon}m^{\varepsilon-1/2}\|\varphi_n(x)\|_{W_2^1}.
\end{equation}
Здесь мы воспользовались интегрированием по частям и учли тот
факт, что $\varphi_n(0)=\varphi_n(\pi)=0$.

Так как из асимптотических формул \eqref{efas} следует, что
$\varphi'_n(x)=n\eta_n(x)+u(x)y_n(x)$, то
\(\|\varphi_n(x)\|_{W_2^1}\le C_{\theta,R}n^{1-\theta}\eta_n\),
где \(\|\{\eta_n\}\|_{l_2}\le C_{\theta, R}\). Значит,
$$\|S_m\varphi_n(x)\|_{C}\le
C_{\theta,R,\varepsilon}m^{\varepsilon-1/2}n^{1-\theta}\eta_n.$$
Итак, первое слагаемое в \eqref{11} не превосходит $$
\|S_m(g_k-g_{N})\|_{C}\le C_{\theta,
R,\varepsilon}\|f\|_{L_2}\left(\suml_{n=1}^{k}m^{2\varepsilon-1}n^{2-2\theta}\eta_n^2\right)^{1/2}\le
C_{\theta,
R,\varepsilon}\|f\|_{L_2}m^{\varepsilon-1/2}k^{1-\theta}.
$$

Шаг 5 завершен.

Перейдем ко второму слагаемому в \eqref{11}. Действие оператора
$S_m$ на функции $g_{N}$ оценивается точно так же, как в
\eqref{12}:
$$
\|S_mg_{N}(x)\|_{C}\le
C_{\varepsilon}m^{\varepsilon-1/2}\|g_{N}(x)\|_{W_2^1}.
$$
В шаге 2 мы доказали оценку: \(\|P_{N}(u)\|_{L_2\to W_2^1}\le
C_{\theta,R}\). Так как \(g_{N}=P_{N}f\), то
\(\|g_{N}\|_{W_2^1}\le C_{\theta, R}\|f\|_{L_2}\). Значит,
\begin{equation}
\label{S1}  \|S_mg_{N}\|_{C}\le C_{\theta,
R,\varepsilon}\|f\|_{L_2}m^{\varepsilon-1/2}.
\end{equation}
Из неравенств \eqref{11}, \eqref{S} и \eqref{S1} теперь сразу
вытекает, что
\begin{equation} \label{13}\|B_mg_k\|_{C}\le C_{\theta,
R,\varepsilon}\|f\|_{L_2}m^{\varepsilon-1/2}k^{1-\theta}.
\end{equation}
 Для завершения доказательства
теоремы осталось положить $M=N^2$ и заметить, что для любого
\(m\ge M\) из \eqref{10} и \eqref{13} следует оценка \eqref{1}
(нужно взять в этих неравенствах $k=[\sqrt{m}]+1$). Теорема
полностью доказана.

Автор благодарит  проф. А.А.Шкаликова и доц. А.М.Савчука за
полезные замечания.

\medskip


\begin{thebibliography}{99}
\addcontentsline{toc}{cont}{{\bf Список литературы}}

\bibitem{Tr} Трибель Х. Теория функциональных
пространств // М.: Мир, 1986.

\bibitem{SSh1} Савчук А.~М., Шкаликов А.~А. Операторы
Штурма--Лиувилля с сингулярными потенциалами // Матем. заметки, Т.
66. \No 6, 1999, С. 897--912.


\bibitem{SSh2} Савчук А.~М., Шкаликов А.~А. Операторы Штурма-Лиувилля
с потенциалами --- распределениями // Труды Московского Мат.
Общества, Т.64, 2003, С. 159--219.

\bibitem{SSh3} Савчук  A.~M., Шкаликов А.~А. О собственных значениях оператора
Штурма--Лиувилля с потенциалами из пространств Соболева.// Матем.
заметки, Т.80, \No 6, 2006, С. 864--884.

\bibitem{Mar} Марченко В.~А. Операторы
Штурма-Лиувилля и их приложения. Киев: Наукова думка, 1977.

\bibitem{Il} Ильин В.~А. О необходимом условии равносходимости с
тригонометрическим рядом спектрального разложения произвольной
суммируемой функции. // Дифф. уравнения, Т. 21, \No 3, 1985, С.
371--379.

\bibitem{VS} Винокуров В.~А., Садовничий В.~А. Равномерная
равносходимость ряда Фурье по собственным функциям первой краевой
задачи и тригонометрического ряда Фурье// Докл. АН, Т. 380, \No 6,
2001, С. 731--735.

\bibitem{R} Гомилко А.~М., Радзиевский Г.~В. Равносходимость рядов
по собственным функциям обыкновенных
функционально-дифференциальных операторов // Докл. АН, Т. 316, \No
2, 1991, С. 265--269.

\bibitem{Sad} Садовничая И.~В. О равносходимости разложений в ряды по тригонометрической системе и
по собственным функциям оператора Штурма--Лиувилля с
потенциалом--распределением // Докл. АН, Т. 392, \No 2, 2003, С.
170--173.

\bibitem{S} Savchuk A.~M. Uniform asymptotic formulae
for eigenfunctions of Sturm--Liouville operators with singular
potentials// arXiv 0801.1950, 2008.

\bibitem{Ka} Като Т. Теория возмущений линейных операторов // М.:
Мир, 1972.

\bibitem{Ba} Бари Н.~К. Тригонометрические ряды.
М.: Государственное издательство физ.-мат. литературы, 1961.



\end{thebibliography}
\end{document}